\documentclass[11pt]{article}

\usepackage{latexsym, amssymb, amsmath, amscd, amsfonts, epsfig, graphicx, colordvi,verbatim,ifpdf}

\newtheorem{thm}{Theorem}[section]

\def\pf{\noindent{\it Proof.} }
\def\qed{\nopagebreak\hfill{\rule{4pt}{7pt}}
\medbreak}

\def\qed{\nopagebreak\hfill{\rule{4pt}{7pt}}
\medbreak}

\title{The third order Ramanujan's mock theta functions}
\author{}

\begin{document}

\begin{center}

{\large \bf Partition Identities for Ramanujan's  Third Order\\[2pt]  Mock Theta Functions}\\[15pt]
\end{center}

\vskip 5mm

\begin{center}
{ William Y. C. Chen}$^{1}$,    {Kathy Q. Ji}$^{2}$,  and {Eric H. Liu}$^3$
\vskip 2mm

   Center for Combinatorics, LPMC-TJKLC\\
   Nankai University, Tianjin 300071, P.R. China

   \vskip 2mm

    $^1$chen@nankai.edu.cn, $^2$ji@nankai.edu.cn, $^3$eric@cfc.nankai.edu.cn
\end{center}

\vskip 6mm \noindent {\bf Abstract.}  We find two involutions on partitions that lead to
partition
identities for Ramanujan's third order mock theta functions
$\phi(-q)$ and $\psi(-q)$. We also give an involution for Fine's partition identity on the mock theta function $f(q)$. The two classical identities of
Ramanujan on third order mock theta functions are
consequences of these partition identities.  Our
combinatorial constructions also apply to Andrews' generalizations
of Ramanujan's identities.

\medskip


\medskip


\section{Introduction}

This paper is concerned with the following three
mock theta functions of order 3 defined by Ramanujan,
\begin{align}
f(q)&=\sum_{n=0}^\infty \frac{q^{n^2}}{(-q;q)_n^2},\\[2pt]
\phi(q)&=\sum_{n=0}^\infty \frac{q^{n^2}}{(-q^2;q^2)_n},\\[2pt]
\psi(q)&=\sum_{n=1}^\infty \frac{q^{n^2}}{(q;q^2)_n}.
\end{align}

Mock theta functions have been extensively studied, see, for example,
  Andrews \cite{Andrews-1989-283}, Fine \cite[Chapters 2-3]{Fine-1988}, Gordon and McIntosh
\cite{Gordon-2010}, and Ono \cite{Ono-2009}. These functions not
only have remarkable analytic properties, but also are  closely
connected to the theory of partitions, see, for example, Agarwal
\cite{Agarwal-2004}, Andrews \cite{Andrews-1966-2,
Andrews-2005-4666}, Andrews and Garvan \cite{ Andrews-1989-242},
Andrews, Eriksson, Petrov, and Romik \cite{Andrews-2007-545}, and
Choi and Kim \cite{Choi-2010-25}.

In this paper, we find two involutions on partitions that imply two
 partition identities for Ramanujan's third order mock theta
functions $\phi(-q)$ and $\psi(-q)$.
We also give an involution for  Fine's partition theorem on the mock theta function $f(q)$. These three partition identities lead to the following two
identities \eqref{r-1} and \eqref{r-2} of Ramanujan
\begin{align}\label{r-1}
&\phi(-q)-2\psi(-q)=f(q),\\
&\phi(-q)+2\psi(-q)=\frac{(q;q)_\infty}{(-q;q)_\infty^2},\label{r-2}
\end{align}
where we have adopted the standard notation
\begin{align}
(a;q)_n&=(1-a)(1-aq)\cdots(1-aq^{n-1}),\\
(a;q)_\infty&=\prod_{n=0}^\infty (1-aq^n),\qquad |q|<1.
\end{align}
 The first proofs of   \eqref{r-1} and \eqref{r-2} were given by Watson   \cite{Watson-1936}.
  Fine \cite[p. 60]{Fine-1988} found another
 proof by using transformation formulas.

\medskip

  Andrews \cite{Andrews-1966-64} defined the following functions as generalizations of Ramanujan's mock
 theta functions and he later found that these generalizations were already  in  Ramanujan's ``lost" notebook \cite{Andrews-1979-89},
\begin{align}
f(\alpha;q)&=\sum_{n=0}^\infty \frac{q^{n^2-n}\alpha^n}{(-q;q)_n(-\alpha;q)_n},\\
\phi(\alpha;q)&=\sum_{n=0}^\infty \frac{q^{n^2}}{(-\alpha q;q^2)_n},\\
\psi(\alpha;q)&=\sum_{n=1}^\infty \frac{q^{n^2}}{(\alpha;q^2)_n}.
\end{align}

When $\alpha=q$, the above functions reduce to Ramanujan's mock
theta functions. Furthermore, Andrews showed these three functions
turn out to be  mock theta functions for $\alpha=q^r$, where $r$ is any positive
integer. More importantly, Ramanujan's
identities \eqref{r-1} and \eqref{r-2} can be extended to
the functions $f(\alpha;q)$, $\phi(\alpha;q)$ and $\psi(\alpha;q)$,
\begin{align}\label{rg-1}
\phi(-\alpha;-q)-(1+\alpha q^{-1})\psi(-\alpha;-q)&=f(\alpha;q),\\
\label{rg-2} \phi(-\alpha;-q)+(1+\alpha
q^{-1})\psi(-\alpha;-q)
&=\frac{(q;q)_\infty}{(-q;q)_\infty(-\alpha;q)_\infty},
\end{align}
see Andrews  \cite[p. 78, Eqs. (3a)--(3b)]{Andrews-1966-64}.
Clearly, the above identities \eqref{rg-1} and \eqref{rg-2}
specialize to \eqref{r-1} and \eqref{r-2}  by setting $\alpha=q$.

The connection between Ramanujan's third order mock theta function
and the theory of partitions was first explored by Fine.  In
\cite[p. 55, Chapter 3]{Fine-1988}, he derived the following
identity from his transformation formula:
\begin{equation}\label{fine-1}
f(q)=1+\sum_{k\geq
1}\frac{(-1)^{k-1}q^k}{(-q;q)_k}=1+\frac{1}{(-q;q)_\infty}\sum_{k\geq
1}(-1)^{k-1}q^k(-q^{k+1};q)_\infty.
\end{equation}
In fact, \eqref{fine-1} can be easily  established from the
combinatorial definition \eqref{Int-f} of $f(q)$.  The following partition identity
  for $f(q)$   can be deduced from \eqref{fine-1}.

 \begin{thm}[Fine]\label{ParIden-1}Let
$p_{do}(n)$  denote the number of partitions of $n$ into distinct
parts with the smallest part being odd. Then
 \begin{equation}\label{pi-1}
(-q;q)_\infty f(q)=1+2\sum_{n\geq 1}p_{do}(n)q^n.
\end{equation}
\end{thm}

We obtain the following two partition identities
from our involutions.

\begin{thm}\label{add-1}We have
\begin{equation} (-q;q)\phi(-q)=1+\sum_{n=1}^\infty
p_{do}(n)q^n+\sum_{k=1}^\infty (-1)^k q^{k^2}.\label{add-1-e}
\end{equation}
 \end{thm}

 \begin{thm}\label{add-2}We have
\begin{equation} 2(-q;q)\psi(-q)=-\sum_{n=1}^\infty
p_{do}(n)q^n+\sum_{k=1}^\infty (-1)^k q^{k^2}.\label{add-2-e}
\end{equation}
 \end{thm}
 Since the generating function for $p_{do}(n)$ is easy to compute, it would be interesting to establish the
above relations as $q$-series identities without resort to
partitions.

It can be seen that the above partition identities lead to
Ramanujan's identities \eqref{r-1} and \eqref{r-2}. It follows
from \eqref{add-1-e} and \eqref{add-2-e} that
\begin{align*}
&(-q;q)_\infty \phi(-q)-2(-q;q)_\infty \psi(-q)\\[2pt]
&\quad=1+\sum_{n=1}^\infty p_{do}(n)q^n+\sum_{k=1}^\infty (-1)^k
q^{k^2}\\
&\quad \quad \quad +\sum_{n=1}^\infty
p_{do}(n)q^n-\sum_{k=1}^\infty (-1)^k q^{k^2}\\
&\quad =1+2\sum_{n\geq 1}p_{do}(n)q^n\\
&\quad =(-q;q)_\infty f(q),
\end{align*}
where the last equality is a consequence of \eqref{pi-1}. So we
obtain the identity \eqref{r-1} by dividing both sides by
$(-q;q)_\infty$.

 In view of
\eqref{add-1-e} and \eqref{add-2-e}, we find
\begin{align*}
&(-q;q)_\infty\phi(-q)+2(-q;q)_\infty \psi(-q)\\
&\quad =1+\sum_{n=1}^\infty p_{do}(n)q^n+\sum_{k=1}^\infty (-1)^k
q^{k^2}\\
&\quad \quad \quad -\sum_{n=1}^\infty
p_{do}(n)q^n+\sum_{k=1}^\infty (-1)^k q^{k^2}\\
&\quad =1+2\sum_{k\geq 1}(-1)^k
q^{k^2}\\
&\quad =\frac{(q;q)_\infty}{(-q;q)_\infty},
\end{align*}
where the last equality follows from  Gauss' identity
\begin{equation}\label{Gauss}
\sum_{k=-\infty }^\infty
(-1)^kq^{k^2}=\frac{(q;q)_\infty}{(-q;q)_\infty}.
\end{equation}
This yields Ramanujan's identity \eqref{r-2} after dividing both sides  by $(-q;q)_\infty$.

In fact, we can deduce two   partition identities for $\phi(-q)$ and $\psi(-q)$ analogous to
Fine's identity for $f(q)$ by employing the following
partition theorem of Bessenrodt and Pak
\cite{Bessenrodt-2004-1139} which extends a theorem of Fine in
\cite[Theorem 5]{Fine-1948}.
It is worth mentioning that there are other involutions which also
imply this partition theorem, see, Berndt, Kim and Yee
\cite{Berndt-2010}, Chen and Liu \cite{Chen}, and Yee
\cite{Yee-2010-a, Yee-2010-b}.

\begin{thm}[Bessenrodt and Pak]\label{thm}Let $p^e_{do}(n)$ {\rm(}$p^o_{do}(n)${\rm)} denote the number of partitions of $n$
into even {\rm(}odd\,{\rm)} distinct parts with the smallest part
being odd. Then
\begin{equation}
\sum_{n=1}^\infty [p^e_{do}(n)-p^o_{do}(n)]q^n=\sum_{k=1}^\infty
(-1)^k q^{k^2}.
\end{equation}
\end{thm}

\medskip

In light of Theorem \ref{thm}, we can restate Theorems
\ref{add-1} and \ref{add-2} as follows.

\begin{thm}\label{ParIden-2}
 We have
\begin{equation}\label{pi-2}
(-q;q)_\infty \phi(-q)=1+2\sum_{n\geq 1}p^e_{do}(n)q^n.
\end{equation}
\end{thm}

\begin{thm}\label{ParIden-3}
We have
\begin{equation}\label{pi-3}
(-q;q)_\infty \psi(-q)=-\sum_{n\geq 1}p^o_{do}(n)q^n.
\end{equation}
\end{thm}

This paper is organized as follows.  In Section 2, we provide an involution for Fine's theorem \ref{ParIden-1}. In Sections
3 and 4, we give two involutions that lead to   partition
identities \eqref{add-1-e} and \eqref{add-2-e} for $\phi(-q)$ and
$\psi(-q)$.  Section 5 is devoted to  partition
identities for $f(\alpha q;q),\phi(-\alpha q;-q)$, and
$\psi(-\alpha q;-q)$ based on our involutions, which imply Andrews' identities \eqref{rg-1} and
\eqref{rg-2}.

\section{An involution for Fine's theorem}

In this section we give an involution for Fine's partition
theorem. Let $P$  denote the set of  partitions, and let $D$
denote the set of partitions with distinct parts.  The
 rank of a partition $\lambda$, denoted by $r(\lambda)$,  is
defined as the largest part minus the number of parts, as
introduced  by Dyson \cite{Dyson-1944}. The empty partition is
assumed to have rank zero. A pair of partitions $(\lambda,\mu)$ is
called a bipartition of $n$ if $|\lambda|+|\mu|=n.$ Fine
\cite[p.49]{Fine-1988} found the following combinatorial
interpretation for $f(q)$:
\begin{equation}\label{Int-f}
f(q)=1+\sum_{\lambda \in P}
(-1)^{r(\lambda)}q^{|\lambda|},
\end{equation}
from which we can construct an involution to prove Theorem
\ref{ParIden-1}.

\medskip

 {\noindent \it Proof of Theorem \ref{ParIden-1}.} For a bipartition $(\lambda,\mu) \in P\times D$ , let $s(\mu)$ denote
the smallest part of $\mu$, $m(\lambda)$ denote the number of
occurrences of the largest part of $\lambda$ and
 $l(\lambda)$ denote the number of (positive) parts of $\lambda$. Let $U$ be the set of
 two classes of bipartitions $(\emptyset,\mu)$ and bipartitions $(\lambda,\mu)$, where
$\lambda=(1,1,\ldots,1)$ and $s(\mu)>l(\lambda)$, and let $V$ be
the set of
 bipartitions $(\lambda,\mu)\in
P\times D$ of $n$ except for  bipartitions in $U$. We shall
construct an involution $\Upsilon$ on the set $V$. The following
two cases are considered.

\begin{itemize}
\item[(1)]If $s(\mu)\leq m(\lambda)$, then delete the smallest
part in $\mu$ and add $1$ to each of the first $s(\mu)$ parts
$\lambda_1,\lambda_2,\ldots, \lambda_{s(\mu)}$ of $\lambda$.
\item[(2)]If $s(\mu)> m(\lambda)$, then subtract
 $1$ from each of the first $m(\lambda)$ parts
$\lambda_1,\lambda_2,\ldots, \lambda_{m(\lambda)}$ of $\lambda$
and add a part of size $m(\lambda)$ to $\mu$.
\end{itemize}

It is easy to check that the above mapping is an involution.
Moreover, $\Upsilon$ changes the parity of the rank of $\lambda$
in $V$.

Let $(\lambda, \mu)$  be a bipartition  in $U$. It is easily seen
that if $l(\lambda)$ is even, then  $r(\lambda)$ is odd. In this
case, we obtain a bipartition $(\emptyset,\nu)$ with $s(\nu)$
being even by moving all the parts of  $\lambda$ to $\mu$ as a
single part which cancels with $(\lambda,\mu)$.  In the case that
$l(\lambda)$ is odd, we see that $r(\lambda)$ is even. Thus we get
bipartition $(\emptyset,\nu)$ with $s(\nu)$ being odd by moving
all the parts of  $\lambda$ to $\mu$ as a single part.  Now, we
are left with two types of bipartitions $(\emptyset,\nu)$ such
that $s(\nu)$ is odd, which correspond to the right  side of
\eqref{pi-1}. This completes the proof. \qed

Here is an example. There are eight bipartitions
$(\lambda,\mu)\in P\times D$ of $4$ with the
rank of $\lambda$ being even,
$$((3),(1)),((2,1),(1)),((1,1,1),(1)),((1),(2,1)),((2,2),\emptyset),
((1),(3)),(\emptyset,(3,1)),(\emptyset,(4)).$$
 Meanwhile, there are
six bipartitions $(\lambda,\mu)\in P\times D$ of $4$ with
the rank of $\lambda$ being odd,
$$((4),\emptyset),((3,1),\emptyset),((2,1,1),\emptyset),((2),(2)),
((1,1),(2)),((1,1,1,1),\emptyset),$$ and there is only one
partition of $4$ into distinct parts with the smallest part being
odd, i.e., $(3,1)$.

\noindent The involution $\Upsilon$ gives the following
correspondence:
\begin{align*}
\begin{array}{ccc}
((3),(1))\leftrightarrow ((4),\emptyset),
&((2,1),(1))\leftrightarrow ((3,1),\emptyset),
&((1,1,1),(1))\leftrightarrow ((2,1,1),\emptyset),\\[6pt]
((1),(2,1))\leftrightarrow
((2),(2)),&((2,2),\emptyset)\leftrightarrow ((1,1),(2)).
\end{array}
\end{align*}
 For the
remaining four bipartitions $((1),(3))$, $(\emptyset,(3,1))$,
$(\emptyset,(4))$, and $((1,1,1,1),\emptyset)$, we can  transform
$((1),(3))$  to $(\emptyset,(3,1))$ and transform
$(\emptyset,(4))$ to $((1,1,1,1),\emptyset)$.

\section{A partition identity  for $\phi(-q)$}

In this section, we shall prove the partition identity for $\phi(-q)$ as stated in Theorem \ref{add-1}. Let us begin with an interpretation
of $\phi(-q)$ given by Fine \cite[p.49]{Fine-1988}.
Let $DO$ denote the set of   partitions with
distinct odd parts. Fine showed that
\begin{equation}\label{Int-phi}
\phi(-q)=1+\sum_{\lambda \in
DO}(-1)^{\frac{\lambda_1+1}{2}}q^{|\lambda|}.
\end{equation}
 Note that Choi and Kim found another interpretation of
 $\phi(q)$ in terms of  $n$-color partitions
 \cite[Theorem
3.1]{Choi-2010-25}.

Hence we have
\begin{equation}\label{add-3}
(-q;q)_\infty \phi(-q)=\sum_{\mu \in
D}q^{|\mu|}+\sum_{(\lambda,\mu) \in DO\times
D}(-1)^{\frac{\lambda_1+1}{2}}q^{|\lambda|+|\mu|}.
\end{equation}

In order to deal with the sum $\sum_{k\geq 1} (-1)^k q^{k^2}$ on
the right hand side of \eqref{add-1-e},  special attention has to
be paid to certain
  bipartitions
$Q_k=(2k-1,2k-3,\ldots,1)$, which  is a partition of $k^2$.

Let $s_o(\mu)$ ($s_e(\mu)$) denote the smallest odd (even) part of
$\mu$, and let $W$ denote the set of bipartitions $(\lambda, \mu)$
of $n$  in $DO\times D$  except for those  of the form $(Q_k,
\emptyset)$ and those bipartitions $(\lambda,\mu)=((1),\mu)$ with
$s_o(\mu)+1<s_e(\mu)$. Obviously, there is a cancellation between
the set of bipartitions  $(\lambda,\mu)=((1),\mu)$ with
$s_o(\mu)+1<s_e(\mu)$ and the set of partitions
 into distinct part with the smallest part
being even  in the first summand of \eqref{add-3}. Consequently,
the remaining partitions for the first summand  in \eqref{add-3}
give the  sum $\sum_{n=1}^\infty p_{do}(n)q^n$ in \eqref{add-1-e}.

Finally, to prove Theorem \ref{add-1}, we are required to
construct an involution on $W$, denoted by $\Phi$, which changes
the parity  of $(\lambda_1+1)/2$.

{\noindent  The involution $\Phi$.} For a partition $\lambda
\in DO$, let $c(\lambda)$ denote the maximum
number of consecutive odd parts of $\lambda$ starting with the first part. For example, let $\lambda=(11,9,7,3)$, then
$c(\lambda)=3.$ The involution $\Phi$ consists of two parts.

\noindent Part I of $\Phi$. If $2c(\lambda)\geq s_e(\mu)$, then
remove the smallest even part from $\mu$, and add $2$ to each of
the first $s_e(\mu)/2$ parts $\lambda_1,\lambda_2,\cdots,
\lambda_{s_e(\mu)/2}$ of $\lambda$.

If $2c(\lambda)< s_e(\mu)$ and $\lambda_{c(\lambda)}>1$, then subtract $2$ from each of the first $c(\lambda)$ parts $\lambda_1,\lambda_2,\ldots,\lambda_{c(\lambda)}$ of $\lambda$, and add a part of size $2c(\lambda)$ to $\mu$.

It is easy to see the above process is well defined and bipartitions which can not be paired by the
 involution are those bipartitions
 $(Q_k, \mu)$ where $s_e(\mu)>2k$. This is the task of the second part of the involution.

\noindent Part II of $\Phi$. If $s_e(\mu)> s_o(\mu)+(2k-1)$, then
delete the smallest odd part of $\mu$ and delete the part $2k-1$
from $Q_k$. Then
 add an even part of size $s_o(\mu)+(2k-1)$  to $\mu$.

If $s_e(\mu)\leq s_o(\mu)+(2k-1)$, then split the smallest even part of $\mu$ into two parts, one is of  size $2k+1$
and the other is of size $s_e(\mu)-(2k+1)$. Observe that $s_e(\mu)-(2k+1)$   is smaller than $s_o(\mu)$. Then add a part of size $2k+1$ to $Q_k$ and add a part of size  $s_e(\mu)-(2k+1)$
to $\mu$.

 So we have obtained an involution $\Phi$. It is readily
seen that this involution changes the parity  of $(\lambda_1+1)/2$.

For example, when $n=7$,  there are $6$ bipartitions
$(\lambda,\mu)\in DO\times D$ with $\lambda_1 \equiv 1\bmod{4}$,
$$((5),(2)),((5,1),(1)),((1),(6)),((1),(4,2)),((1),(3,2,1)),((1),(5,1)).$$
On the other hand,  there are $5$ bipartitions $(\lambda,\mu)\in
DO\times D$ with $\lambda_1 \equiv 3\bmod{4}$, that is,
$$((7),\emptyset),((3,1),(2,1)),((3,1),(3)),((3),(4)),((3),(3,1)),$$
and there is only one partition of $7$ into distinct parts with
the smallest part being even, namely, $(5,2)$.

 The involution $\Phi$ gives the following pairs
 of bipartitions:
\begin{align*}
\begin{array}{ccc}
((5),(2))\leftrightarrow ((7),\emptyset),
&((5,1),(1))\leftrightarrow ((3,1),(2,1)),
 &((1),(6))\leftrightarrow ((3,1),(3)),\\[6pt]
((1),(4,2))\leftrightarrow ((3),(4)),
&((1),(3,2,1))\leftrightarrow ((3),(3,1)).
\end{array}
\end{align*}
For the remaining bipartition $((1),(5,1))$, we can construct a partition into
distinct part with smallest part being even, that is, $(5,2)$.

\section{A partition identity for $\psi(-q)$}

The aim of this section is to prove the partition identity
\eqref{add-2-e} for $\psi(-q)$. There is also a combinatorial
interpretation of $\psi(-q)$ given by Fine \cite[p.49]{Fine-1988}.
Let $OC$ denote the set of partitions of $n$ into odd parts
without gaps. Fine \cite[p.49]{Fine-1988} showed that
\begin{equation}\label{Int-psi}
\psi(-q)=\sum_{\lambda \in
OC}(-1)^{l(\lambda)}q^{|\lambda|}.
\end{equation}
Note that Agarwal \cite{Agarwal-2004, Agarwal-2007} found two  combinatorial interpretations for $\psi(q)$ by using $q$-difference equations.

Let $D^0$ denote the set of partitions with distinct parts where the zero
part is allowed. So
 the number of partitions of $n$ in the
set $D^0$ is twice  the number of partitions of $n$ in the set
$D$. Let $R$ denote the set of bipartitions $(\lambda,\mu)\in
OC\times D^0$ except for those  of the form $(Q_k, \emptyset)$ and
those bipartitions $(\lambda,\mu)=((1),\mu)$ with
$s_e(\mu)+1<s_o(\mu)$. It is clear that
 these excluded bipartitions correspond to the right hand side of  \eqref{add-2-e}.
In order to prove \eqref{add-2-e}, it   suffices to construct an involution $\Psi$ on the set $R$.

{\noindent  The involution $\Psi$.}
Let $r_p(\lambda)$ be the largest part of $\lambda$ which occurs
at least twice in $\lambda$, where we let $r_p(\lambda)=\infty$ if
$\lambda$ has no repeated parts. The involution $\Psi$ consists of three parts.

\noindent Part I of $\Psi$. If $s_o(\mu)> r_p(\lambda)$, then
delete one part of size $r_p(\lambda)$ from $\lambda$ and add it
as a part to $\mu$. On the other hand, if $s_o(\mu)\leq
r_p(\lambda)$, then move a part of size $s_o(\mu)$ from $\mu$ to
$\lambda$.

The above process is well defined and the bipartitions not covered by this case  are those bipartitions $(Q_k, \mu)$ for which $s_o(\mu)>2k-1.$
 We  continue to describe the second part
 of $\Psi$.

 \noindent
 Part II of $\Psi$.
 Assume that $(Q_k,\mu)$ is a bipartition such that $s_o(\mu)>2k-1$.
If $\mu$ has a zero part, then we get a bipartition
$(Q_{k-1},\mu^*)$ where $\mu^*$ is obtained from $\mu$ by adding a
part of size $2k-1$ and deleting the zero part.

If $\mu$ has no zero part and  $s_o(\mu)=2k+1$, then we get a
bipartition $(Q_{k+1},\mu^*)$ where $\mu^*$ is obtained from $\mu$
by removing a part of size $2k+1$ and adding a zero part.

It can be seen that the above mapping  is well defined except for  those bipartitions $(Q_k, \mu)$ such that $\mu$ has
no zero part and $s_o(\mu)>2k+1$.
Indeed, it is the object of the third part of
$\Psi$ to deal with these remaining bipartitions.

\noindent Part III of $\Psi$.
If $s_o(\mu)> s_e(\mu)+(2k-1)$, then delete the smallest even part of $\mu$ and delete the part $2k-1$ from $Q_k$,
 and add an odd part of size $s_e(\mu)+(2k-1)$  to $\mu$.

If $s_o(\mu)\leq s_e(\mu)+(2k-1)$, then split the
smallest odd part of $\mu$ into two parts, one is of size $2k+1$
and the other is of size $s_o(\mu)-(2k+1)$, which is less than
$s_e(\mu)$,  add a part of size $2k+1$ to $Q_k$ and add a
part of size  $s_o(\mu)-(2k+1)$ to $\mu$.

It is routine  to check that the map $\Psi$ is
an involution and it changes the parity of the length of $\lambda$.

For  example, when $n=4$  there are six
bipartitions $(\lambda,\mu)\in OC\times D^0$
such that $l(\lambda)$ is odd,
$$((1),(3)),((1,1,1),(1)),((1,1,1),(1,0)),((1),(2,1)),((1),(2,1,0)),((1),(3,0)).$$
In the other case,  there are  five bipartitions $(\lambda,\mu)\in
OC\times D^0$ such that $l(\lambda)$ is even,
$$((3,1),(0)),((1,1,1,1),\emptyset),((1,1,1,1),(0)),((1,1),(2)),((1,1),(2,0)),$$
and there is one partition  of $4$ with distinct parts such that the
smallest part is odd, i.e., $(3,1)$.

The involution  $\Psi$ is illustrated below:   \begin{align*}
\begin{array}{lll}
((1),(3))\leftrightarrow ((3,1),(0)),
&((1,1,1),(1))\leftrightarrow ((1,1,1,1),\emptyset),
 \ ((1),(2,1))\leftrightarrow ((1,1),(2))\\[6pt]
((1),(2,1,0))\leftrightarrow ((1,1),(2,0)),&
((1,1,1),(1,0))\leftrightarrow ((1,1,1,1),(0)).
\end{array}
\end{align*}
For the remaining bipartition $((1),(3,0))$, we can form a partition into
distinct parts with smallest part being odd, that is, $(3,1)$.

As  another example, the involution $\Psi$ also gives the
following correspondence:
$$((9,7,5,3,1),(16,15,8,6,2))\leftrightarrow ((7,5,3,1),(16,15,11,8,6)).$$

\section{Andrews' generalizations}

This section is devoted to proofs of Andrews' identities
\eqref{rg-1} and \eqref{rg-2}. First, we give combinatorial
interpretations for $f(\alpha q;q), \phi(-\alpha q;-q),$ and
$\psi(-\alpha q;-q)$ by extending the arguments of Fine. More
precisely, we have the following partition theoretic
interpretations.

\begin{thm}\label{int-g}We have
\begin{align}\label{Int-f-g}
f(\alpha q;q)&=1+\sum_{\lambda \in P} (-1)^{r(\lambda)}
 \alpha^{\lambda_1}q^{|\lambda|},\\ \label{Int-phi-g}
 \phi(-\alpha q;-q)&=1+\sum_{\lambda \in
DO}
(-1)^{\frac{\lambda_1+1}{2}}\alpha^{\frac{\lambda_1+1}{2}-l(\lambda)}q^{|\lambda|},\\
\label{Int-psi-g}
 \psi(-\alpha q;-q)&=\sum_{\lambda \in
OC}(-1)^{l(\lambda)} \alpha^{l(\lambda)-\frac{\lambda_1+1}{2}}
q^{|\lambda|}.
\end{align}
\end{thm}

\pf Recall that
$$f(\alpha q;q)=\sum_{n=0}^\infty \frac{q^{n^2}\alpha^n}{(-q;q)_n(-\alpha q;q)_n},$$
it is easy to check that
\eqref{Int-f-g} follows from the
Durfee square dissection of a
partition $\lambda \in P$, see Figure 1.

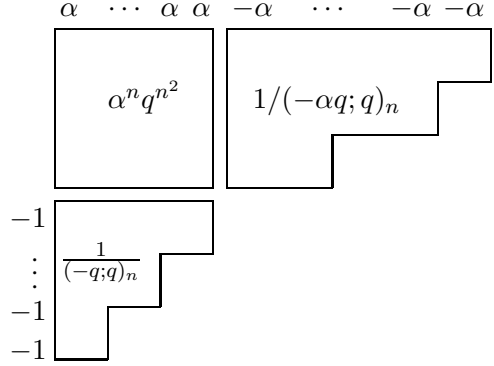
\begin{figure}[h]
\begin{center}
\begin{picture}(220,160)
\put(20,140){\line(1,0){60}} \put(20,140){\line(0,-1){60}}
\put(20,80){\line(1,0){60}} \put(80,140){\line(0,-1){60}}
\put(85,140){\line(1,0){100}} \put(185,140){\line(0,-1){20}}
\put(185,120){\line(-1,0){20}} \put(165,120){\line(0,-1){20}}
\put(165,100){\line(-1,0){40}} \put(125,100){\line(0,-1){20}}
\put(125,80){\line(-1,0){40}} \put(85,80){\line(0,1){60}}
\put(20,75){\line(0,-1){60}} \put(20,75){\line(1,0){60}}
\put(80,75){\line(0,-1){20}} \put(80,55){\line(-1,0){20}}
\put(60,55){\line(0,-1){20}} \put(60,35){\line(-1,0){20}}
\put(40,35){\line(0,-1){20}} \put(40,15){\line(-1,0){20}}
\put(40,110){$\alpha^n q^{n^2}$} \put(22,145){$\alpha$}
\put(72,145){$\alpha$} \put(60,145){$\alpha$}
\put(40,145){$\cdots$}

\put(87,145){$-\alpha$} \put(117,145){$\cdots$}
\put(147,145){$-\alpha$} \put(167,145){$-\alpha$} \put(3,15){$-1$}
\put(3,65){$-1$} \put(3,30){$-1$} \put(11,40){$\cdot$}
\put(11,45){$\cdot$} \put(11,50){$\cdot$} \put(95,110){$1/(-\alpha
q;q)_n$} \put(22,50){$\frac{1}{(-q;q)_n}$}
\end{picture}
\caption{The Durfee square dissection.}
\end{center}
\end{figure}

From the definition of $\phi(\alpha;q)$, we see that
$$\phi(-\alpha q;-q)=\sum_{n=0}^\infty
\frac{(-1)^nq^{n^2}}{(-\alpha q^2;q^2)_n}.$$ The term
$(-1)^nq^{n^2}$ corresponds  to a partition $\pi$ of the form
$Q_n=(2n-1,2n-3,\ldots,3,1)$, which has weight
$(-1)^{(\pi_1+1)/2}$. Moreover,  $1/(-\alpha q^2;q^2)_n$
 is the generating function for partitions $\sigma$
with at most $n$ even parts and with no odd parts. The weight of such a partition is endowed with  weight $(-\alpha)^{\sigma_1/2}$.

 Define
 $\lambda=\pi+\sigma=(\pi_1+\sigma_1, \pi_2+\sigma_2, \ldots)$. We see that $\lambda \in DO$, namely, $\lambda$ is a partition into distinct odd parts. Now,  the weight of $\lambda$ equals $(-1)^{(\lambda_1+1)/2}\alpha^{(\lambda_1
+1)/2-l(\lambda)}$. So
 \eqref{Int-phi-g} has been verified.

For the combinatorial interpretation for $\psi(-\alpha q;-q)$, we
note that
$$\psi(-\alpha q;-q)=\sum_{n=1}^\infty \frac{(-1)^nq^{n^2}}{(-\alpha q;q^2)_n}.$$
 The summand can be expanded as follows
\begin{align*}
\frac{(-1)^nq^{n^2}}{(-\alpha q;q^2)_n}&=\frac{-q}{1+\alpha
q}\frac{-q^3}{1+\alpha q^3}\cdots \frac{-q^{2n-1}}{1+\alpha
q^{2n-1}}\\
&=(-q+\alpha q^{1+1}-\alpha^2 q^{1+1+1}\cdots)(-q^3+\alpha
q^{3+3}-\alpha^2 q^{3+3+3}+\cdots)\\
&\qquad \cdots (-q^{2n-1}+\alpha q^{2(2n-1)}-\alpha^2
q^{3(2n-1)}\cdots).
\end{align*}
 It follows that the summand $(-1)^nq^{n^2}/(-\alpha q;q^2)_n$ is the generating function of
partitions $\lambda$ in $OC$ with the largest part not exceeding
$2n-1$ and with weight $
(-1)^{l(\lambda)}\alpha^{l(\lambda)-(\lambda_1+1)/2}$. This proves
(\ref{Int-psi-g}). \qed

We can extend Fine's partition theorem to Andrews' function
$f(\alpha q;q)$. It can be  seen that the involution $\Upsilon$ in
Section 2 preserves the quantity $\lambda_1+l(\mu)$. Therefore,
from \eqref{Int-f-g} we deduce the following partition theorem.

 \begin{thm}\label{ParIden-1-g}Let $P_{do}$ denote the set of
 partitions into distinct parts with the smallest part being odd.
 Then
 \begin{equation}\label{pi-1-g}
 (-\alpha q;q)_\infty f(\alpha q;q)=1+2\sum_{\nu \in P_{do}}\alpha^{l(\nu)}q^{|\nu|}.
 \end{equation}
\end{thm}

Next, we give a  generalization of Theorem \ref{add-1} to
$\phi(-\alpha q;-q)$. By the combinatorial interpretation
\eqref{Int-phi-g}, we find that
\begin{align}
(-\alpha q;q)_\infty\phi(-\alpha q;-q)&=(-\alpha q;q)_\infty+\sum_{(\lambda,\mu)
\in DO\times D}(-1)^{\frac{\lambda_1+1}{2}}\alpha^{\frac{\lambda_1+1}{2}-l(\lambda)
+l(\mu)}q^{|\lambda|+|\mu|}.
\end{align}
Moreover, we observe that the involution $\Phi$ in Section 3
preserves the quantity of
\[ l(\mu)-l(\lambda)+\frac{\lambda_1+1}{2}.\] Hence we have
\begin{equation}\label{temp-3}
\sum_{(\lambda,\mu) \in DO\times D}(-1)^{\frac{\lambda_1+1}{2}}
\alpha^{\frac{\lambda_1+1}{2}-l(\lambda)+l(\mu)}q^{|\lambda|+|\mu|}
=-\sum_{\nu \in P_{de}}\alpha^{l(\nu)}q^{|\nu|}+\sum_{k=1}^\infty (-1)^k q^{k^2},
\end{equation}
where $P_{de}$ denotes the set of partitions into distinct parts
with the smallest part being even. On the other hand,
\begin{equation}\label{temp-4}
(-\alpha q;q)_\infty=1+\sum_{\nu \in P_{de}}\alpha^{l(\nu)}q^{|\nu|}+\sum_{\nu \in P_{do}}\alpha^{l(\nu)}q^{|\nu|},
\end{equation}
Therefore, from \eqref{temp-3} and \eqref{temp-4} we deduce the following  partition identity for
$\phi(-\alpha q;-q)$.

\begin{thm}\label{ParIden-2-g}We have
\begin{equation}\label{pi-2-g}
(-\alpha q;q)_\infty\phi(-\alpha q;-q)=1+\sum_{\nu \in P_{do}}\alpha^{l(\nu)}q^{|\nu|}+\sum_{k=1}^\infty (-1)^k q^{k^2}.
\end{equation}
\end{thm}

We now proceed to give a generalization of Theorem \ref{add-2} to
$\psi(-\alpha q;-q)$. By the combinatorial interpretation
\eqref{Int-psi-g}, we obtain that
\begin{equation}
(-\alpha;q)_\infty\psi(-\alpha q;-q)=\sum_{(\lambda,\mu) \in
DO\times D^0} (-1)^{l(\lambda)}\alpha^{l(\lambda)-
\frac{\lambda_1+1}{2}+l(\mu)}q^{|\lambda|+|\mu|}.
\end{equation}

It can be verified that the involution $\Psi$ in Section 4 preserves the quantity of
\[ l(\lambda)-\frac{\lambda_1+1}{2}+l(\mu)\]
 and it changes the parity of $l(\lambda)$.
So we get the following partition theorem.

\begin{thm}\label{ParIden-3-g}We have
\begin{equation}\label{pi-3-g}
(-\alpha;q)_\infty\psi(-\alpha q;-q)=-\sum_{\nu \in P_{do}}\alpha^{l(\nu)}q^{|\nu|}+\sum_{k=1}^\infty (-1)^k q^{k^2}.
\end{equation}
\end{thm}

Based on the above partition theorems for $f(\alpha
q;q),\,\phi(-\alpha q;-q)$ and $\psi(-\alpha q;-q)$, we can deduce
Andrews' generalizations of Ramanujan's identities. More
precisely, it follows from  \eqref{pi-2-g} and \eqref{pi-3-g} that
\begin{align*}
&(-\alpha q;q)_\infty\phi(-\alpha q;-q)-(-\alpha;q)_\infty\psi(-\alpha q;-q)\\[2pt]
&\quad =1+\sum_{\nu \in P_{do}}\alpha^{l(\nu)}q^{|\nu|}+\sum_{k=1}^\infty (-1)^k q^{k^2}\\
&\quad \qquad +\sum_{\nu \in P_{do}}\alpha^{l(\nu)}q^{|\nu|}-\sum_{k=1}^\infty (-1)^k q^{k^2}\\
&\quad =1+2\sum_{\nu \in P_{do}}\alpha^{l(\nu)}q^{|\nu|}\\
&\quad = (-\alpha q;q)_\infty f(\alpha q;q),
\end{align*}
where the last equality is a consequence of   identity
\eqref{pi-1-g}. Dividing both sides by $(-\alpha q ;q)_\infty$ yields
\[\phi(-\alpha q;-q)-(1+\alpha)\psi(-\alpha q;-q)=f(\alpha q;q).
\]
Hence we deduce the identity \eqref{rg-1} by  replacing $\alpha
q$  by $\alpha$.

According to \eqref{pi-2-g} and \eqref{pi-3-g}, we have
\begin{align*}
&(-\alpha q;q)_\infty\phi(\alpha q;-q)+(-\alpha;q)_\infty\psi(\alpha q;-q)\\[2pt]
&\quad =1+\sum_{\nu \in P_{do}}\alpha^{l(\nu)}q^{|\nu|}+\sum_{k=1}^\infty (-1)^k q^{k^2}\\
&\qquad \quad -\sum_{\nu \in P_{do}}\alpha^{l(\nu)}q^{|\nu|}+\sum_{k=1}^\infty (-1)^k q^{k^2}\\
&\quad =1+2\sum_{k=1}^\infty (-1)^k q^{k^2}\\
&\quad =\frac{(q;q)_\infty}{(-q;q)_\infty},
\end{align*}
where the last equality follows from Gauss' identity
\eqref{Gauss}. Dividing both sides by $(-\alpha q ;q)_\infty$, we
obtain
\[\phi(-\alpha q;-q)+(1+\alpha)\psi(-\alpha q;-q)=\frac{(q;q)_\infty}{(-q;q)_\infty(-\alpha q;q)_\infty},
\]
we arrive at the identity \eqref{rg-2} by replacing
$\alpha q$ by  $\alpha$.

To conclude, we note that our approach can be
 viewed as combinatorial proofs of  Andrews' identities in the forms multiplied by the factor $(-\alpha q ;q)_\infty$.

\noindent{\bf Acknowledgments.} This work was supported by the 973
Project, the PCSIRT Project and the Doctoral Program Fund of the Ministry of Education, the
National Science Foundation of China.


\begin{thebibliography}{99} \small
\setlength{\itemsep}{-.8mm}

\bibitem{Agarwal-2004}A.K. Agarwal, $n$-Color partition theoretic interpretations of
some mock theta functions, Electron. J. Combin. 11 (2004) N14.

\bibitem{Agarwal-2007}A.K. Agarwal, New combinatorial interpretations
of some mock theta functions, Online J. Analytic Combin. 2 (2007), \#5.


 \bibitem{Andrews-1966-64}G.E. Andrews, On basic
hypergeometric series, mock-theta functions, and partitions (I),
Quart. J. Math. (2) 17 (1966) 64--80.

\bibitem{Andrews-1966-2}G.E. Andrews, On basic
hypergeometric series, mock-theta functions, and partitions (II),
Quart. J. Math. (2) 17 (1966) 132--143.

\bibitem{Andrews-1976} G.E. Andrews, The Theory of Partitions,
Addison-Wesley Publishing Co., 1976.

\bibitem{Andrews-1979-89}G.E. Andrews, An introduction to Ramanujan's ``Lost"
Notebook, Amer. Math. Monthly 86 (1979) 89--108.


\bibitem{Andrews-1989-242}G.E. Andrews and F.G. Garvan, Ramanujan's ``lost" notebook VI: The mock theta conjectures,
Adv. Math. 73 (1989) 242--255.

\bibitem{Andrews-1989-283}G.E. Andrews, Mock theta functions, In Theta Functions, Bowdoin 1987, Part 2, Proc. Symp.
Pure Math., vol. 49, American Mathematical Society, Providence,
RI, (1989) 283--298.



\bibitem{Andrews-2005-4666}G.E. Andrews, Partitions with short
sequences and mock theta functions, Proc. Natl. Acad. Sci. USA 102
(2005) 4666--4671.

\bibitem{Andrews-2007-545}G.E. Andrews, H. Eriksson, F. Petrov and
D. Romik, Integrals, partitions and MacMahon's theorem, J. Combin.
Theory Ser. A 114 (2007) 545--554.

\bibitem{Berndt-2010}B.C. Berndt, B. Kim and A.J. Yee, Ramanujan's Lost Notebook:
Combinatorial proofs of identities associated with Heine's
transformation or partial theta functions, J. Combin. Theory Ser.
A., to appear.



\bibitem{Bessenrodt-2004-1139}B.C. Bessenrodt and I. Pak, Partition congruences by involutions,
European J. Combin. 25 (2004) 1139--1149.

\bibitem{Chen} W.Y.C. Chen and E.H. Liu, A Franklin type involution for squares,   Adv. Appl. Math., to appear.


\bibitem{Choi-2010-25}Y.-S. Choi and B. Kim, Partition identities from third and sixth order mock theta functions, preprint.

\bibitem{Dyson-1944}F.J. Dyson, Some guesses in the theory of
partitions, Eureka {\rm(}Cambridge{\rm)}  8 (1944) 10--15.

\bibitem{Fine-1948}N.J. Fine, Some new results on partitions, Proc. Nat. Acad. Sci. USA
34 (1948) 616--618.

\bibitem{Fine-1988}N.J. Fine,
Basic Hypergeometric Series and Applications, Math. Surveys
27, AMS Providence, 1988.



\bibitem{Gordon-2010}B. Gordon and R. McIntosh, A survey of classsical mock
theta functions, preprint.


\bibitem{Ono-2009}K. Ono, Unearthing the visions of a master: Harmonic
Maass forms and number theory, In: Proceeding of the 2008 Harvard-MIT Current Developments in Mathematics Conference, International Press, Somerville, MA, 2009, pp. 347--454.


\bibitem{Watson-1936} G.N. Watson, The final problem: An account of the mock theta functions, J. London
Math. Soc. 11 (1936) 55--80.



    \bibitem{Yee-2010-a}A.J. Yee, Bijective proofs of a theorem of Fine and related partition identities,
Internat. J. Number Theory, to appear.

\bibitem{Yee-2010-b}A.J. Yee, Ramanujan's partial theta series and parity in partitions,
Ramanujan J., to appear.



\end{thebibliography}
\end{document}